\documentclass[11pt, two side]{article}
\usepackage{amsmath,amssymb}
\usepackage{array}
\usepackage{latexsym}
\usepackage{amsfonts}
\usepackage{graphicx}
\usepackage {hyperref}
\usepackage{indentfirst}
\numberwithin{equation}{section}

\begin{document}
\def \b{\Box}

\begin{center}
{\Large {\bf An elementary method for the fidel codification of
texts written in Romanian language\\[0.2cm]
 }}
{\Large {\bf (O metod\u a elementar\u a  pentru codificarea
fidel\u a a textelor scrise {\^\i}n limba rom\^an\u a)\\[0.2cm]
 }}
\end{center}

\begin{center}
{\bf Mihai IVAN}
\end{center}

\setcounter{page}{1}

\pagestyle{myheadings}

{\bf Abstract.} In the paper \cite{miha}, the Polybius square is
extended at the dimension $7\times 7$ (called, {\bf $7\times
7-$Square}). Choosing of the characters of the $7\times 7-$Square
is determined by the need for fidel codification of  texts written
in Romanian language.

In this paper an elementary method (cryptographic system) for the
codification and decodification  of texts written in Romanian
language is presented.This method is based on the alphabet ${\bf
A}$ and the $7\times 7-$Square.

 The theme presented is at the interference  between
mathematics, Romanian language and cryptography.

Learning mathematics in high school pursues awareness of nature of
mathematics as one dynamic discipline closely related to real
life, by its relevance in everyday life and its role in the
sciences, and in technology.

This paper may be used in compiling the list of optional content
in  the discipline of mathematics (IX class), more precisely at
the subject: "{\it Bijective functions. Application: $7\times
7-$Square and the fidel codification of texts written in Romanian
language}".

Through the organization of learning activities of a set of
lessons with the above subject, the teacher provides solutions for
effective teaching-
learning process, contributing to:\\
$-~$ acquisition  of skills arising from the perspective
mathematics applications in the real-world;\\
$-~$ skills information, research and systematization of
information and personal development of students through
communication, cooperation and  active-learning.{}\\

 {\bf Rezumat}. \^{I}n lucrarea \cite{miha}, p\u atratul lui Polybius este extins la dimensiunea
 $ 7\times 7 $(numit, P\u atratul $ 7\times 7 $). Alegerea caracterelor P\u
atratului  $ 7\times 7 $ este determinat\u a de necesitatea de
codificare fidel\u a a textelor scrise \^{i}n limba rom\^{a}n\u a.

\^{I}n aceast\u a lucrare vom prezenta o metod\u a elementar\u a
(sistem de criptare)  pentru codificarea \c si decodificarea
textelor scrise \^{i}n limba rom\^{a}n\u a. Aceast\u a metod\u a
se bazeaz\u a pe alfabetul ${\bf A}$  \c si P\u atratul $ 7\times
7 .$

 Tema prezentat\u a se afl\u a la
interferen\c ta dintre matematic\u a, limba rom\^{a}n\u a \c si
criptografie.

\^{I}nv\u a\c tarea matematicii \^{i}n liceu urm\u are\c ste con\c
stientizarea naturii matematicii ca o disciplin\u a dinamic\u a
\^{i}n str\^{a}ns\u a leg\u atur\u a cu via\c ta real\u a prin
relevan\c ta sa \^{i}n cotidian \c si prin rolul s\u au \^{i}n \c
stiin\c te \c si  \^{i}n tehnologii.

 Acest articol poate fi
utilizat \^{i}n alc\u atuirea listei con\c tinuturilor unui op\c
tional la disciplina "matematic\u a" (clasa a $IX$-a), mai precis
la tema: "{\it Func\c tii bijective. Aplica\c tie: P\u atratul $
7\times 7$  \c si secretizarea fidel\u a a textelor
 scrise \^{i}n limba rom\^{a}n\u a}".

Prin organizarea unor activit\u a\c ti de \^{i}nv\u a\c tare a
unui set de lec\c tii cu tema de mai sus, profesorul ofer\u a
solu\c tii pentru un proces de predare-\^{i}nv\u a\c tare
eficient, care contribuie la:\\
$-~$ dobâ\^{a}ndirea unor competen\c te  ce decurg din perspectiva
aplica\c tiilor matematicii \^{i}n lumea real\u a;\\
$-~$ dezvoltarea abilit\u a\c tilor de informare, cercetare \c si
sistematizare a informa\c tiilor \c si dezvoltarea personal\u a a
elevilor prin comunicare, cooperare \c si \^{i}nv\u a\c tare
activ\u a.

\section{Introducere}
\indent Rezolvarea unor probleme practice din via\c ta cotidian\u
a, \^{i}n cadrul lec\c tiilor de matematic\u a, joac\u a un rol
determinant \^{i}n \^{i}nv\u a\c tarea unor concepte matematice.
Utilizarea aplica\c tiilor practice ale matematicii \^{i}n diverse
domenii contribuie la formarea deprinderilor de munc\u a
intelectual\u a \^{i}n \^{i}nv\u a\c tare.

\markboth{Mihai Ivan}{ O metod\u a elementar\u a pentru
codificarea fidel\u a a textelor scrise {\^\i}n limba rom\^an\u a}

Proiectarea \c si transpunerea \^{i}n practic\u a a unor unit\u
a\c ti de \^{i}nv\u a\c tare care s\u a propun\u a elevilor
situa\c tii de \^{i}nva\c tare prin intermediul aplica\c tiilor
practice vor conduce la dob\^{a}ndirea unor cuno\c stin\c te
calitativ superioare, precum \c si la cre\c sterea performan\c
telor la matematic\u a.

Criptografia este \c stiin\c ta comunic\u arii informa\c tiei sub
o form\u a securizat\u a (\cite{atan, pirc}). Odat\u a cu
transmiterea de mesaje, oamenii au avut nevoie ca acestea ajunse
\^{i}n posesia unor persoane indezirabile s\u a nu poat\u a fi
decodificate. Dezvoltarea tot mai rapid\u a a telecomunica\c
tiilor \c si a comunica\c tiilor electronice au impus tot mai mult
crearea de noi tehnici criptografice. Tehnicile criptografice sunt
folosite pentru a securiza comunica\c tiile derulate prin
intermediul re\c telelor de calculatoare, pentru a efectua pl\u
ati on-line \c si a implementa scheme de vot electronic, \^{i}n
cadrul telefoanelor mobile sau cardurilor bancare etc. Ast\u azi,
un num\u ar mare de persoane utilizeaz\u a algoritmi \c si tehnici
din domeniul criptografiei, care au la baz\u a o serie de concepte
\c si metode matematice.

P\u atratul $ 7\times 7 $ este asociat alfabetului ${\bf A}$
(\cite{miha}). Alfabetul ${\bf A}$ con\c tine $31$ litere mici ale
alfabetului limbii rom\^{a}ne \c si $ 18 $ caractere care permit
codificarea fidel\u a a textelor scrise \^{i}n limba rom\^{a}n\u a
(de exemplu, fraze, strofe din poezii, citate \c si cuget\u ari
celebre etc.). Cifrarea oric\u arui caracter se realizeaz\u a
aleg\^{a}nd num\u arul de dou\u a cifre care corespunde a\c sez\u
arii caracterului \^{i}n p\u atrat. Cifrarea unei litere mari se
face prin scrierea al\u aturat\u a a dou\u a numere de dou\u a
cifre \c si anume: primul num\u ar este asociat caracterului
${\cal L}$ (adic\u a num\u arul $ 54 $ care indic\u a ac\c tiunea
de transformare a literi mici \^{i}n liter\u a mare), iar al
doilea num\u ar este cel care corespunde literei mici \^{i}n p\u
atrat.

{\bf Exemplu}. Mesajul \^{i}n clar: "{\it Rom\^{a}nia}"  se
transform\u a dup\u a cifrarea bazat\u a pe P\u atratul $ 7\times
7$  \^{i}n mesajul cifrat: " $ 54~34~31~26~51~27~22~11$".

Con\c tinutul \c stiin\c tific al lucr\u arii este accesibil \c si
u\c sor de \^{i}n\c teles pentru elevii din \^{i}nv\u a\c t\u
am\^{a}ntul preuniversitar (clasele a IX-a$-$XII-a).

Proiectarea \c si organizarea unei activit\u a\c ti didactice
bazat\u a pe elementele de con\c tinut ale acestei lucr\u ari este
benefic\u a din urm\u atoarele motive:\\
  $-~~$se prezint\u a \^{i}ntr-un mod elementar utilitatea conceptelor matematice (de exemplu: coresponden\c te \^{i}ntre cuvinte \c si secven\c te de numere,
  func\c tie  bijectiv\u a);\\
 $-~~$contribuie la familiarizarea elevilor cu unele aplica\c tii simple ale criptografiei. Prin descifrarea mesajului trimis de c\u atre destinatar, elevii au
 satisfac\c tia reu\c sitei;\\
  $-~~$prin dob\^{a}ndirea abilit\u a\c tilor de cifrare, elevii vor \^{i}nv\u a\c ta ceva util, interesant \c si atractiv.

\section{No\c tiuni de criptografie}

$\indent$ Primele texte cifrate descoperite p\^an\u a {\^\i}n
prezent dateaz\u a de circa $4000$ de ani \c si provin din Egiptul
antic. {\^I}n Grecia antic\u a, scrierile cifrate erau folosite
{\^\i}nc\u a din secolul al V-lea a.Hr.
Polybius{\footnote{Polybius ($\approx$ 207-120 a.Hr.), istoric
grec}} a creat un tabel de cifrare {\^\i}n form\u a de p\u atrat
de dimensiune $5\times 5$ (cunoscut sub numele de {\it p\u atratul
lui Polybius}). {\^I}n Roma antic\u a secretul informa\c tiilor
politice \c si militare se realiza folosind scrierea secret\u a
(de exemplu, {\it cifrul lui Cezar}{\footnote{Julius
Caesar($\approx$ 100-44 a.Hr.), lider politic \c si general
roman}}).

{\^I}n criptografia clasic\u a (numit\u a, {\it criptografia
pre-computa\c tional\u a}) au fost create diverse sisteme de
cifrare bazate pe P\u atratul lui Polybius care sunt utilizate \c
si ast\u azi (\cite{pirc}).

Termenul de {\it criptografie}  {\^\i}nseamn\u a {\it scriere
secret\u a}. Cuv\^antul criptografie este format din cuvintele
grece\c sti {\bf cryptos} ({\it ascuns}) \c si {\bf grafie} ({\it
scriere}). {\^I}n sens restr\^ans, termenul de criptografie
desemneaz\u a numai opera\c tia de cifrare \c si descifrare
legal\u a.

{\^I}n cele ce urmeaz\u a prezent\u am unele no\c tiuni de baz\u a
din domeniul criptografiei (\cite{atan, pope}).

{\bf Criptografia} ({\it cryptography}) este \c stiin\c ta cre\u
arii \c si p\u astr\u arii mesajelor secrete astfel {\^\i}nc\^at
descifrarea lor s\u a fie imposibil\u a de c\u atre persoane
neautorizate.

\indent  Un {\bf mesaj} ${\cal M}$ (numit {\it text {\^\i}n clar})
scris {\^\i}n limba rom\^an\u a este mesajul ce urmeaz\u a a fi
secretizat. Un {\bf mesaj cifrat} ${\cal C}$ (numit {\it
criptogram\u a})  este mesajul secretizat, care este accesibil
numai persoanelor autorizate.

{\bf Criptare/cifrare} ({\it encryption})$~{\cal E}$ este
procedeul de "ascundere" a unui mesaj {\^\i}n clar {\^\i}n mesajul
secretizat. Avem: $~~~{\cal E}({\cal M})= {\cal C}.$

Procedeul de secretizare a unui text este realizat de c\u atre
expeditor.

 {\bf Decriptare/descifrare} ({\it decryption})$~{\cal D}$ este procedeul de
reg\u asire a mesajului {\^\i}n clar din mesajul cifrat. Avem:
  $~~~ {\cal D}({\cal C})= {\cal D}({\cal E}({\cal M}))={\cal M}.$

Desecretizarea unui text criptat este realizat\u a de c\u atre
destinatar.

Prin {\bf algoritm criptografic/cifru} se desemneaz\u a o mul\c
time de transform\u ari inversabile prin care mul\c timea
mesajelor se transform\u a {\^\i}n mul\c timea mesajelor cifrate
\c si invers. Cifrul  este construit cu ajutorul a dou\u a func\c
tii (func\c tia de criptare $ {\cal E}$ \c si func\c tia de
decriptare $ {\cal D}$).

{\bf Cheia criptografic\u a} ({\it key}) $~{\cal K}$ este m\u
arimea (secret\u a) necesar\u a realiz\u arii cript\u arii \c si
decript\u arii. Cheia de cifrare este reprezentat\u a printr-un
num\u ar, cuv\^ant, \c sir numeric etc. \c si care reglementeaz\u
a opera\c tia de cifrare. Algoritmul care realizeaz\u a opera\c
tiile de criptare \c si decriptare se nume\c ste {\it sistem de
criptare}.

Pentru aplicarea unui sistem de criptare se procedeaz\u a astfel:\\[0.1cm]
  $-~$un mesaj {\^\i}n forma sa original\u a   este un  text {\^\i}n
  clar;\\
  $-~$expeditorul rescrie mesajul ${\cal M},$  folosind un procedeu (algoritm) cunoscut numai de el (\c si de destinatar). Se spune c\u a el cripteaz\u a (cifreaz\u a)
  mesajul ${\cal M},$  ob\c tin\^and un text cifrat ${\cal C}$;\\
  $-~$destinatarul prime\c ste textul cifrat ${\cal C}$  \c si {\^\i}l decripteaz\u a, cunosc\^and algoritmul folosit pentru
  criptare;\\
  $-~$se presupune c\u a expeditorul \c si destinatarul stabilesc de comun acord, {\^\i}ntr-o faz\u a preliminar\u a, detaliile procedurii de criptare \c si de
  decriptare.

\section{P\u atratul  $7\times 7 $}

Consider\u am un alfabet ${\bf A}={\bf A}_{1}\cup {\bf A}_{2}$
format din 49 caractere, unde:\\
$\bullet~~~{\bf A}_{1}=\{~a, b, c, d, e, f, g, h, i, j, k, l, m,
n, o, p, q, r, s, t, u, v, w, x, y, z,$  \u a,  \^ i,  \^ a,  \c s,  \c t $\}$\\
 este {\it mul\c{t}imea format\u a din 31 litere} (litere mici ale alfabetului limbii rom\^ ane);\\
   $\bullet~~~{\bf A}_{2} $    este {\it mul\c timea format\u a din 18  caractere} (simboluri  sau semne grafice de punctua\c tie) care au urm\u atoarele
   semnifica\c tii:\\
    ${\cal L}~~~~~\ldots~~~~~$        {\it litera mic\u a care urmeaz\u a devine} {\bf liter\u a
    mare};\\
    ${\sqcup}~~~~~\ldots~~~~~$  un {\bf spa\c tiu liber} (delimitator de spa\c tiu);\\
    $/~~~~~\ldots~~~~~$   {\bf r\^and nou};\\
    $,~~~~~\ldots~~~~~$  {\bf virgul\u a};\\
    {-}$ ~~~~~\ldots~~~~~$  {\bf cratim\u a};\\
    $? ~~~~~\ldots~~~~~$      {\bf semnul {\^\i}ntreb\u arii};\\
    $! ~~~~~\ldots~~~~~$     {\bf semnul exclam\u arii};\\
$" ~~~~~\ldots~~~~~$   {\bf  ghilimele pentru scrierea unui text};\\
$" ~~~~~\ldots~~~~~$ {\bf ghilimele pentru {\^\i}ncheierea unui
text};\\
$; ~~~~~\ldots~~~~~$     {\bf punct \c si virgul\u a};\\
$- ~~~~~\ldots~~~~~$   {\bf linia de pauz\u a} sau {\bf linia de
dialog};\\
$\centerdot ~~~~~\ldots~~~~~$      {\bf punct};\\
$: ~~~~~\ldots~~~~~$             {\bf dou\u a puncte};\\
$^{,} ~~~~~\ldots~~~~~$             {\bf  apostrof};\\
$( ~~~~~\ldots~~~~~$              {\bf parantez\u a deschis\u
a};\\
$)~~~~~\ldots~~~~~$              {\bf parantez\u a {\^\i}nchis\u
a};\\
$\& ~~~~~\ldots~~~~~$   {\bf \c si};\\
$@~~~~~\ldots~~~~~$   {\bf simbol pentru e-mail}.\\

Cele  $49$ de litere \c{s}i simboluri ale alfabetului {\bf A} sunt
a\c{s}ezate \^{i}ntr-un tabel {\^\i}n form\u a de p\u atrat,
format $7$ linii \c{s}i $7$ coloane ($L_{i}$ \c{s}i $C_{j}$ pentru
$i,j=\overline{1,7}$), numit "{\bf P\u atratul} $7\times 7$" sau "
{\bf P\u atratul lui Polybius extins}" (vezi tabelul de mai jos).

\[
\quad\begin{array}{|r|c|c|c|c|c|c|c|} \hline
  & C_{1} & C_{2} & C_{3} & C_{4} & C_{5} & C_{6} & C_{7}\\ \hline

L_{1}  & a & b & c & d & e & f & g  \cr \hline

L_{2}  & h & i & j & k & l & m & n \cr \hline

L_{3}  & o & p & q & r & s & t & u  \cr \hline

L_{4}  & v & w & x & y & z & $\u{a}$  & $\^{i}$ \cr \hline

L_{5} & $\^{a}$ & $\c{s}$ & $\c{t}$  & {\cal L} & \sqcup & / & ,
\cr \hline

L_{6} & $-$ & ? & ! & " & " & ; & - \cr \hline

L_{7} & \centerdot & : & ^{,} & ( & ) & $\&$ & @  \cr \hline
\end{array}.
\]

{\bf Observa\c tie.} Fiecare din simbolurile $\cal{L}, \sqcup, / $
ne indic\u a executarea unei ac\c tiuni asupra unei litere mici,
cuv\^{a}nt sau r\^{a}nd (mai precis, transformarea unei litere
mici a alfabetului ${\bf A}_{1} $ \^{i}n litera mare corespunz\u
atoare, l\u{a}sarea unui spa\c tiu liber dup\u a un cuv\^{a}nt
scris sau trecerea la un r\^{a}nd nou dup\u a un r\^{a}nd
scris).\hfill$\Box$

 Orice caracter din P\u atratul $7\times 7$  este
"identificat" cu o pereche de numere corespunz\u atoare pozi\c
tiei caracterului \^{i}n p\u atrat (mai precis, perechea
respectiv\u a este format\u a din num\u arul liniei \c si num\u
arul coloanei pe care se afl\u a litera sau simbolul). De exemplu,
avem coresponden\c tele urm\u atoare:\\[0.1cm]
$p~~\leftrightarrow~~(3,2);~~~z~~\leftrightarrow~~(4,5);~~~M~~\leftrightarrow~~(5,4);~~~\sqcup~~\leftrightarrow~~(5,5);~~~f~~\leftrightarrow~~(1,6).$

    \^{I}n acest mod, alfabetului  ${\bf A}$ i se asociaz\u a un nou alfabet notat cu  {\bf A(P)} \c si numit {\it alfabetul perechilor  asociat  P\u atratului $7\times 7$}.
    Avem:\\[-0.4cm]
$${\bf A(P)}= \{~(i,j)~|~i,j=\overline{1,7}\}, $$
unde $(i,j)~$ este perechea  asociat\u a caracterului din linia
$L_{i}$ \c si coloana  $C_{j}$ a P\u atratului $7\times 7$.

\section{Utilizarea numerelor \^{i}n codificarea fidel\u a  \c
si decodificarea textelor scrise \^{i}n limba rom\^{a}n\u{a}}

{\^I}n acest paragraf prezent\u am un sistem de criptare pentru
codificarea \c si decodificarea textelor, care se bazeaz\u a pe
alfabetul $~{\bf A}$ \c si P\u atratul $ 7 \times 7.$ Acesta este
un cifru de substitu\c tie (substitution cipher). Cifrarea fiec\u
arui caracter al textului \^{i}n clar ($ M $) este substituit cu
un num\u ar de dou\u a cifre sau grup de numere de c\^{a}te dou\u
a cifre \^{i}n textul cifrat ($ C $), descifrarea
realiz\^{a}ndu-se prin aplicarea substitu\c tiei inverse.

{\bf Pentru a secretiza un text scris {\^\i}n limba rom\^ an\u a} se vor respecta urm\u atoarele conven\c tii:\\
 $-~~~${\it textul {\^\i}n clar este format din cuvinte, propozi\c tii sau fraze} care con\c tin litere mici \c si litere
mari ale alfabetului limbii rom\^ane, semne de punctua\c tie \c
si spa\c tii libere {\^\i}ntre cuvinte;\\
$-~~~${\it textul este interpretat ca un singur cuv\^ant}
(inclusiv  spa\c tiul liber este un caracter distinct).

\begin{center}
\textbf{Codificarea textelor, utiliz\^{a}nd P\u atratul  $7\times
7$}
\end{center}

Dac\u a num\u arul natural $ i\in\{p, p+1, \ldots, q-1, q\}\subset
{\bf N},$ atunci vom scrie uneori $i=\overline{p,q}.$

Opera\c tia de codificare  a unui text, utiliz\^{a}nd alfabetul
${\bf A}$ \c si P\u atratul $7\times 7$, const\u a
\^{i}n:\\
 $\bullet~~${\it fiec\u arei perechi din alfabetul $\bf{A}(\bf P)$  i se asociaz\u a num\u arul $~\overline{ij}~$  format din
 cifrele $ i,j=\overline{1,7}$};\\
 $\bullet~~${\it fiec\u arei litere mari  care corespunde literei mici din pozi\c tia $(m,n)$ i se asociaz\u a succesiunea numerelor $54$  \c{s}i $~\overline{mn},$
  unde  $ m=\overline{1,4}$ \c si  $ n=\overline{1,7} $  sau $ m=5$   \c si $ n=\overline{1,3}.$}

De exemplu, avem coresponden\c tele urm\u atoare:\\
$a~\leftrightarrow~(1,1)~\leftrightarrow~11;~~~~~n~\leftrightarrow~(2,7)~\leftrightarrow~27;~~~~w~\leftrightarrow~(4,2)~\leftrightarrow~42;~~~~~$\c
s$~~\leftrightarrow~(5,2)~
\leftrightarrow~52;$\\[0.1cm]
$M\leftrightarrow~(5,4)~\leftrightarrow~54;~~~~~~\sqcup~\leftrightarrow~(5,5)~\leftrightarrow~55;~~~~~~C~\leftrightarrow~54~13;~~~~~U~\leftrightarrow~
54~37;~~~~~$\c{T}$~\leftrightarrow~54~53.$

Codificarea unui mesaj ${\cal M}$  se exprim\u a matematic cu
ajutorul unei func\c tii bijective $ f$  definit\u a pe alfabetul
$\bf{A}(\bf P)$ format din $ 49 $ perechi asociate P\u atratului
$7\times 7$ cu valori \^{i}n mul\c timea  ${\bf B}=\{\overline{ij}
~|~ i,j=\overline{1,7}\} $ a numerelor formate din dou\u a cifre $
i $ \c si $ j .$ Mai precis:\\[-0.4cm]
$$f: {\bf A}({\bf P})~\longrightarrow~{\bf
B},~~~~~~~f(i,j)=\overline{ij},$$
 unde $ (i,j) $ este perechea asociat\u a
caracterului din linia $L_{i}$  \c si coloana $L_{j}$ a P\u
atratului $7\times 7.$\\

 Dac\u a "identific\u am" alfabetul  $\bf{A}(\bf P)$ cu alfabetul ${\bf A},$ atunci func\c tia $ f$ (numit\u a {\it func\c tie de codificare})  se poate descrie cu ajutorul tabelului
urm\u ator (numit, {\it tabelul asociat func\c tiei $ f $}):
\[
\quad\begin{array}{|r|c|c|c|c|c|c|c|c|c|c|c|c|c|c|c|c|c|c|c|c|}
\hline
 {\bf A} & a & b & c & d & e & f & g &  h & i & j & k & l & m & n & o & p & q & r & s \\ \hline
{\bf B}  & 11 & 12 & 13 & 14 & 15 & 16 & 17 &  21 & 22 & 23 & 24 &
25 & 26 & 27 & 31 & 32 & 33 & 34 & 35  \cr \hline
\end{array}
\]
\[
\quad\begin{array}{|r|c|c|c|c|c|c|c|c|c|c|c|c|c|c|c|c|c|c|c|}
\hline
 {\bf A} & t & u & v & w & x & y & z & $\u{a}$  & $\^{i}$ & $\^{a}$ & $\c{s}$ & $\c{t}$  & {\cal L} & \sqcup & / & ,& $-$ & ? & ! \\ \hline
{\bf B}  & 36 & 37 & 41 & 42 & 43 & 44 & 45 & 46 & 47 & 51 & 52 &
53 & 54 & 55 & 56 & 57 &  61 & 62 & 63 \cr \hline
\end{array}
\]
\[
\quad\begin{array}{|r|c|c|c|c|c|c|c|c|c|c|c|} \hline
 {\bf A} & " & " &  ; & - & \centerdot & : & ^{,} & ( & ) & $\&$ & @ \\ \hline
{\bf B}  & 64 & 65 & 66 & 67 & 71 & 72 & 73 & 74 & 75 & 76 & 77
\cr \hline
\end{array}.
\]

Prin  {\it cuv\^{a}nt de lungime} $n,$  notat $c=c_{1}c_{2}\ldots
c_{n}$ vom \^{i}n\c telege o succesiune format\u a din  $ n $
caractere $ c_{1}, c_{2},~\ldots~, c_{n},$ unde $
c_{i},~i=\overline{1,n}$ este un element al alfabetului ${\bf A} $
(except\^{a}nd simbolul $\cal{L}$ ) sau o liter\u a mare a
alfabetului limbii rom\^{a}ne.

Func\c tia  $ f $ se extinde la mul\c timea cuvintelor de lungime
finit\u a, astfel:\\
 {\it dac\u a $ c= c_{1}c_{2}\ldots c_{n} $ este un cuv\^{a}nt de
lungime $n,$ atunci $ f(c)$ este o succesiune de numere de dou\u a
cifre definit\u a prin $ f(c)= f(c_{1})f(c_{2})\ldots f(c_{n}), $
unde  $ f(c_{i}) $ este un num\u ar de dou\u a cifre sau o
succesiune format\u a din dou\u a numere de c\^{a}te dou\u a
cifre}. Deci:

$c= c_{1}c_{2}\ldots c_{n}~~~\longrightarrow~~~f(c)=
f(c_{1})f(c_{2})\ldots f(c_{n}).$

Observ\u am c\u a $ f(c)$ este o secven\c t\u a de forma $
\overline{i_{1}j_{1}}~ \overline{i_{2}j_{2}}\ldots
\overline{i_{n+q}j_{n+q}},$ unde $ q (0\leq q \leq n )$
reprezint\u a num\u arul literelor mari care intr\u a \^{i}n
componen\c ta cuv\^{a}ntului $c=c_{1}c_{2}\ldots c_{n}.$

De exemplu, cuv\^{a}ntul:\\
  $-~~~c^{\prime}=$"{\it U.V.T.}" este format din $ 6 $  caractere, iar secven\c ta  $ f(c^{\prime}) $ este format\u a din $ 9 $
  numere;\\
  $-~~~c^{\prime\prime}=$"{\it Diana \c{s}i Ana}" este format din $ 12$  caractere, iar secventa $ f(c^{\prime\prime}) $  este format\u a din  $ 14 $
  numere.

 {\bf Pentru a aplica metoda bazat\u a pe alfabetul ${\bf A}$, P\u atratul $ 7\times 7 $  \c si utilizarea numerelor de dou\u a cifre  \^{i}n codificarea textelor},
 se vor aplica urm\u atoarele reguli generale:\\
 $-~~${\it fiec\u arui caracter} (care intr\u a \^{i}n componen\c ta unui cuv\^ ant) {\it i se asociaz\u a prin
intermediul unei func\c tii bijective un num\u ar format din dou\u a cifre}. Mai precis:\\
 $\bullet~~$unei litere mici sau simbol (care reprezint\u a un semn de punctua\c tie sau o ac\c tiune) i se asociaz\u a un num\u ar format din dou\u a
 cifre;\\
 $\bullet~~$unei litere mari i se asociaz\u a o succesiune format\u a din dou\u a numere de c\^ ate dou\u a cifre (al doilea num\u ar este cel asociat  literei mici
 corespunz\u atoare).\\
 $-~~${\it textul codificat este format dintr-o succesiune de numere de dou\u a
 cifre}.

{\bf Exemplul 1.} Cuv\^{a}ntul  $ c= $"{\it teorema lui Pitagora}"
este format din $ 20 $ caractere, iar secven\c ta $ f(c)$  este
format\u a din $ 21 $ numere. Codificarea cuv\^{a}ntului $ c $
prin func\c tia  $ f $ este:\\[0.1cm]
{\it teorema lui
Pitagora}$~~~\longrightarrow~~~36~15~31~34~15~26~11~55~25~37~22~55~54~32~22~36~11~17~31~34~11.$

{\bf Exemplul 2.} Codific\u am urm\u atorul citat al profesorului
{\it Dan Barbilian}{\footnote{Dan Barbilian (1895-1961), poet \c
si matematician rom\^{a}n.}}(cunoscut ca poet sub numele de {\it
Ion Barbu}):"

{\it Exist\u a undeva, \^{i}n domeniul \^{i}nalt al geometriei, un
loc luminos unde se \^{i}nt\^{a}lne\c ste cu poezia}.".\\[0.1cm]
 $-~~$Textul \^{i}n clar este format din  $ 92 $
 caractere.\\
$-~~$Pentru litera  $E$ vom folosi secven\c ta de numere $ 54~15;$
pentru litera  \u a  folosim num\u arul $46;$  pentru caracterul
$,$ (virgul\u a) folosim num\u arul  $ 57$ etc.\\
 $-~~$Textul codificat este format din urm\u atoarele $ 93 $  numere:
 "\\[0.1cm]
$54~15~43~22~35~36~46~55~37~27~14~15~41~11~57~55~47~27~55~14~31~26~15~27~22~37~25~55~47~27~11~25~36~55$\\[0.1cm]
$11~25~55~17~15~31~26~15~36~34~22~15~22~57~55~37~27~55~25~31~13~55~25~37~26~22~27~31~35~55~37~27~14~15$\\[0.1cm]
$55~35~15~55~47~27~36~51~25~27~15~52~36~15~55~13~37~55~32~31~15~45~22~11~71$".

{\bf Exemplul 3.} Codific\u am urm\u atorul mesaj:\\
"{\it  Teorema reciproc\u a a teoremei lui Thales}".\\[0.1cm]
$-~~$Se procedeaz\u a \^{i}n mod analog ca la Exemplul 2.\\
 $-~~$Textul este compus din $ n=39 $
 caractere.\\
 $-~~$Textul codificat
este format din urm\u atoarele $ 41 $  numere:"\\
 $54~36~15~31~34~15~26~11~55~34~15~13~22~32~34~31~13~46~55~11~55~36~15~31~34~15~26~15~22~55~25$\\
$~~~~~~~~~37~22~55~54~36~21~11~25~15~35$".

\begin{center}
\textbf{Decodificarea textelor, utiliz\^{a}nd P\u atratul
$7\times 7$}
\end{center}

Opera\c tia de decodificare  a unui text, utiliz\^{a}nd alfabetul
${\bf A}$ \c si P\u
atratul  $7\times 7$ const\u{a} \^{i}n:\\
 $\bullet~~${\it fiec\u arui num\u ar $\overline{ij}$  format din
 cifrele $i,j=\overline{1,7},$  i se
asociaz\u a perechea $~(i,j)\in \bf{A}(\bf P)$ (iar aceasta  se
identific\u a cu caracterul corespunz\u
ator din P\u atratul  $7\times 7$)};\\
 $\bullet~~${\it fiec\u arei succesiuni de tipul $~ 54~\overline{mn}, $ i se asociaz\u a litera mare care corespunde literei mici din pozi\c tia $(m,n)$ a P\u atratului
  $7\times 7,$ unde  $ m=\overline{1,4}$ \c si  $ n=\overline{1,7}~$  sau $~ m=5$   \c si $ n=\overline{1,3}.$}

De exemplu, avem coresponden\c tele urm\u atoare:\\
$13~\leftrightarrow~(1,3)~\leftrightarrow~c;~~~~~22~\leftrightarrow~(2,2)~\leftrightarrow~i;~~~~~46~\leftrightarrow~(4,6)~\leftrightarrow~$\u
a$;~~~~~62~\leftrightarrow~(6,2)~
\leftrightarrow~?\\[0.1cm]
74~\leftrightarrow~(7,4)~\leftrightarrow~);~~~~~~36~\leftrightarrow~(3,6)~\leftrightarrow~t;~~~~~~54~45~\leftrightarrow~Z;~~~~~54~41~\leftrightarrow~
V;~~~~54~47~\leftrightarrow~$\^{I}.

Decodificarea unui mesaj ${\cal C}$ se exprim\u a matematic cu
ajutorul unei func\c tii bijective $ g $ definit\u a pe  $ {\bf
B}=\{\overline{ij} ~|~ i,j=\overline{1,7}\}~$ cu valori \^{i}n
 $\bf{A}(\bf P).$  Mai precis:\\[-0.4cm]
$$g: {\bf
B}~\longrightarrow~ {\bf A}({\bf P}),~~~~~~~g(\overline{ij})=
(i,j),$$
 unde $ (i,j) $ este perechea asociat\u a
caracterului din linia $L_{i}$  \c si coloana $L_{j}$ a P\u
atratului $7\times 7.$

Dac\u a "{\it identific\u am}" alfabetul $ \bf{A}(\bf P)$ cu
alfabetul ${\bf A},$ atunci {\it func\c tia $ g $ (numit\u a {\it
func\c tie de decodificare}) se poate descrie cu ajutorul unui
tabel care se ob\c tine prin inversarea liniilor tabelului asociat
func\c tiei $f.$}

Prin  {\it secven\c t\u a admis\u a de lungime $n (n\geq 2), $}
notat\u a $ \alpha = \overline{i_{1}j_{1}}~
\overline{i_{2}j_{2}}\ldots \overline{i_{n}j_{n}},$ vom \^{i}n\c
telege o succesiune format\u a din numere de forma
$~\overline{i_{1}j_{1}}, \overline{i_{2}j_{2}}, \ldots,
\overline{i_{n}j_{n}} $ (unde,  $ i_{k}, j_{k}=\overline{1,7} $ \c
si $ k=\overline{1,n} $), care \^{i}ndepline\c ste urm\u atoarele
condi\c tii:\\[0.1cm]
$(i)~~~\overline{i_{k}j_{k}}\neq 54 $ (adic\u a ultimul num\u ar
al secven\c tei nu poate fi egal cu  $ 54 $ );\\[0.1cm]
$(ii)~~$ {\it dac\u a $~\overline{i_{k}j_{k}} = 54, $ atunci $
~\overline{i_{k+1}j_{k+1}}~(1\leq k \leq n-1 ) $} {\it este
acceptat pentru $~i_{k+1}=\overline{1,4} $ \c si
$~j_{k+1}=\overline{1,7}~$ sau $~i_{k+1}=5 $ \c si
$~j_{k+1}=\overline{1,3}.$}

Func\c tia $ g$  se extinde la mul\c timea secven\c telor admise
de lungime finit\u a, astfel:\\
{\it dac\u a $ \alpha = \overline{i_{1}j_{1}}~
\overline{i_{2}j_{2}}\ldots \overline{i_{n}j_{n}} $ este o
secven\c t\u a admis\u a de lungime $ n, $  atunci $ g(\alpha) $
este o succesiune de perechi de numere, definit\u a prin $~
g(\overline{i_{1}j_{1}}) g(\overline{i_{2}j_{2}}) \ldots
g(\overline{i_{n}j_{n}}),$ unde\\[0.1cm]
 $g(\overline{i_{m}j_{m}})=(i_{m}, j_{m} ),~ 1\leq  m\leq  n.~$} Deci:\\[-0.4cm]
$$\alpha = \overline{i_{1}j_{1}} \overline{i_{1}j_{1}} \ldots
\overline{i_{n}j_{n}} ~~~\longrightarrow~~~g(\alpha)=
g(\overline{i_{1}j_{1}})g(\overline{i_{1}j_{1}}) \ldots
g(\overline{i_{n}j_{n}}).$$

Observ\u am c\u a secven\c ta admis\u a  $ \alpha=
\overline{i_{1}j_{1}}  \overline{i_{2}j_{2}} \ldots
\overline{i_{n}j_{n}} $  poate con\c tine  $q$   numere egale cu $
54 $ ( $ q $ este egal cu num\u arul de apari\c tii ale ac\c
tiunii ${\cal L}$). Utiliz\^{a}nd func\c tia bijectiv\u a dintre
alfabetul ${\bf A}({\bf P}) $ \c si alfabetul $ {\bf A},~
g(\alpha) $ se interpreteaz\u a ca un cuv\^{a}nt de lungime $
n-q.$

De exemplu, secven\c ta:\\
 $\alpha^{\prime}=$ "$54~15~37~13~25~22~14$ " este admis\u a de lungime $7$ , iar cuv\^{a}ntul $ g(\alpha) $  are  $ 6 $
 caractere.\\
 $\alpha^{\prime}=$ "$11~54~62~27~11$" nu este admis\u a, deoarece secven\c ta  $ 54~ 62 $   nu este acceptat\u
 a.

{\bf Exemplul 4.}  (a) Secven\c ta  $\alpha=$ "$
35~31~34~22~27~77~15~61~37~41~36~71~34~31$ " este admis\u a \c si
este format\u a din $ 14 $ numere, iar cuv\^{a}ntul $ g(\alpha) $
este format din $ 14 $ caractere. Decodificarea secven\c tei $
\alpha $  prin func\c tia de decodificare $ g $ este:
"$ sorin@e$-$uvt.ro".$\\
(b)  Fie mesajul criptat $ {\cal C},~$  reprezentat prin secven\c
ta format\u a din $ 19 $ numere:\\
 " ${\cal C}: ~~~ 12 ~ 22~
27~ 31~ 26 ~ 37~ 25 ~ 55~ 25~ 37~ 22 ~ 55~ 54~ 27~ 15 ~ 42~ 36 ~31
~ 27$".

Pentru aplicarea func\c tiei de decodificare $ g $  putem folosi
tabelul urm\u ator:
\[
\quad\begin{array}{|r|c|c|c|c|c|c|c|c|c|c|c|c|c|c|c|c|c|c|c|}
\hline
 {\cal C} & 12 & 22 & 27 & 31 & 26 & 37 & 25 & 55 & 25 & 37 & 22 & 55 & 54 & 27 & 15 & 42 & 36 & 31 & 27 \\ \hline
{\cal M}  & b & i & n & o & m & u & l &  & l  & u & i &  & & N & e
& w & t & o & n \cr \hline
\end{array}
\]

Ob\c tinem un cuv\^{a}nt $ c $  format din $18$ caractere,
unde:$~~~ c= $" {\it binomul lui Newton}".

{\bf Exemplul 5.} Decodific\u am urm\u atorul mesaj criptat,
reprezentat prin urm\u atoarea secven\c t\u a format\u a din $ 52
$ de numere:"\\[0.1cm]
$54~27~~37~55~15~35~36~15~55~14~15~35~36~37~25~55~35~46~55~16~22~22~55~12~37~27~71~55~54~36~34$\\
$~~~~~~~~~15~12~37~22~15~55~35~46~55~16~22~22~55~25~11~55~13~15~41~11~71$".
Prin decodificare se descoper\u a urm\u atorul citat al poetului
{\it Tudor Arghezi} {\footnote{Tudor Arghezi (1880-1967), poet \c
si prozator
rom\^{a}n.}}:"\\[0.1cm]
{\it Nu este destul s\u a fii bun. Trebuie s\u a fii bun la
ceva}".

{\bf Exemplul 6.} Decodific\u am urm\u atorul mesaj criptat,
reprezentat prin urm\u atoarea secven\c t\u a format\u a din $ 134
$ de numere:"\\[0.1cm]
$54~32~15~55~13~51~27~14~55~27~15~61~27~41~51~25~36~31~34~15~11~26~55~13~37~55~23~31~13~37~25~57~56~54$\\[0.1cm]
$52~31~32~36~15~11~22~72~55~64~54~47~27~13~15~34~13~37~22~15~61~26~22~55~26~22~23~25~31~13~37~25~63~65$\\[0.1cm]
$56~54~52~22~61~11~13~37~26~55~53~22~61~11~52~55~47~26~32~25~22~27~22~55~14~31~34~22~27~53~11~57~56~54$\\[0.1cm]
$14~11~34~55~27~37~61~53~22~55~13~37~32~34~22~27~14~55~13~22~34~13~37~26~16~15~34~22~27~53~11~63$".

Prin decodificare descoperim o strof\u a scris\u a de Ion Barbu:"

$~~~~~~~~~~~~~~~~~~~~~~~~~~~~~~~~~~~$ {\it Pe c\^{a}nd
ne-nv\^{a}ltoream cu jocul},\\
 $~~~~~~~~~~~~~~~~~~~~~~~~~~~~~~~~~~~~~~~~$ {\it \c{S}opteai:"\^{I}ncercuie-mi mijlocul!"}\\
$~~~~~~~~~~~~~~~~~~~~~~~~~~~~~~~~~~~~~~~~$ {\it \c{S}i-acum
\c{t}i-a\c{s} \^{i}mplini dorin\c{t}a,}\\
$~~~~~~~~~~~~~~~~~~~~~~~~~~~~~~~~~~~~~~~~$ {\it Dar nu-\c{t}i
cuprind circumferin\c{t}a !}".

{\bf Observa\c tie.} Prin alegerea unei alte aranj\u ari a
caracterelor \^{i}n P\u atratul $7\times 7$ se va ob\c tine un nou
cifru de substitu\c tie.\hfill$\Box$\\

{\bf Concluzii.} Sistemul  de cifrare prezentat î\^{i}n aceast\u a
lucrare permite redarea exact\u a \^{i}n textele criptate a
valorilor stilistice ale frazelor scrise \^{i}n limba rom\^{a}n\u
a.

Secretizarea fidel\u a a textelor scrise \^{i}n limba rom\^{a}n\u
a se realizeaz\u a utiliz\^{a}nd alfabetul ${\bf A},$ P\u atratul
$ 7\times 7 $  \c si secven\c tele de numere de dou\u a cifre din
mul\c timea ${\bf B}=\{\overline{ij} ~|~ i,j=\overline{1,7}\}$.

Unul din obiectivele principale ale \^{i}nv\u a\c t\u arii
matematicii \^{i}n liceu este dezvoltarea abilit\u a\c tilor de
\^{i}n\c telegere \c si cunoa\c stere  a conceptelor matematice
pentru a le folosi cu succes \^{i}n situa\c tii concrete.

Pentru \^{i}ndeplinirea acestui deziderat, profesorul poate
propune cursul optional (la nivel de disciplin\u a) cu titlul:
 "{\it \^{I}nv\u a\c t\u am matematic\u a prin
aplica\c tiile ei \^{i}n via\c ta cotidian\u a}".

Con\c tinutul acestui articol poate fi utilizat ca  suport de curs
\^{i}n conceperea \c si realizarea unei unit\u a\c ti de \^{i}nv\u
a\c tare ca parte component\u a a cursului op\c tional men\c
tionat mai sus \c si care se adreseaz\u a elevilor din clasa a
IX-a. Activit\u a\c tile de predare-\^{i}nv\u a\c tare-evaluare
organizate pe parcursul desf\u a\c sur\u arii cursului pot fi
realizate, folosind metode active \c si interactive (observarea
dirijat\u a, problematizarea, jocul didactic, referatul, proiectul
etc.).

\vspace*{0.2cm}

Author's adresses\\[-0.2cm]

Mihai Ivan\\
West University of Timi\c soara,\\
Departamentul pentru Preg{\u a}tirea Personalului Didactic (DPPD),\\
 4, Bd. V. P{\^a}rvan, 300223, Timi\c soara, Romania\\
E-mail: mihai.ivan@e-uvt.ro\\

 \end{document}